\documentclass[11pt]{amsart}
 
\usepackage{amsfonts}
\usepackage{amssymb}
\usepackage{amscd}
\input xy
\xyoption{all}
 
\newcommand{\marginlabel}[1]%
  {\mbox{}\marginpar{\raggedleft\hspace{0pt}\bfseries\sf#1}}

\def\QQ{{\mathbb Q}}

\def\OO{{\mathcal O}}

\def\F{\mathcal{F}}

\def\I{\mathcal{I}} 
\def\J{\mathcal{J}}

\def\Pic0{{\rm Pic}^0(X)}

\theoremstyle{plain}
\newtheorem*{introtheorem}{Theorem}
\newtheorem*{introcorollary}{Corollary}
\newtheorem{theorem}{Theorem}[section]
\newtheorem{proposition}[theorem]{Proposition}
\newtheorem{corollary}[theorem]{Corollary}
\newtheorem{lemma}[theorem]{Lemma}

\theoremstyle{definition}
\newtheorem{definition}[theorem]{Definition}
\newtheorem{remark}[theorem]{Remark}
\newtheorem{example}[theorem]{Example}

\newtheorem{conjecture/question}[theorem]{Conjecture/Question}
  
\newtheorem{remark/definition}[theorem]{Remark/Definition} 
 
\theoremstyle{remark}

\begin{document}
 
\title{Regularity on abelian varieties I}
 
\author[G. Pareschi]{Giuseppe Pareschi} 
\address{Dipartamento di Matematica, Universit\`a di Roma, Tor Vergata, V.le della 
Ricerca Scientifica, I-00133 Roma, Italy}
\email{{\tt pareschi@mat.uniroma2.it}}

\author[M. Popa]{Mihnea Popa$^1$}
\footnotetext[1]{The second author was partially 
supported by a Clay Mathematics Institute
Liftoff Fellowship during the preparation of this paper.}
\address{Department of Mathematics, Harvard University,
One Oxford Street, Cambridge, MA 02138, USA }
\email{{\tt mpopa@math.harvard.edu}}

\maketitle


\markboth{G. PARESCHI and M. POPA}
{\bf REGULARITY ON ABELIAN VARIETIES I}

\section{\bf Introduction}

This is the first in a series of papers meant to introduce a notion of regularity
on abelian varieties and more general irregular varieties. This notion, called 
\emph{Mukai regularity}, is based on Mukai's concept of Fourier transform, and in a very 
particular form (called \emph{Theta regularity}) it parallels and strengthens the usual 
Castelnuovo-Mumford regularity with respect to polarizations on abelian varieties.
Mukai regularity turns out to have a large number of applications, ranging 
from basic properties of linear series on abelian varieties and defining equations 
for the ideals of their subvarieties, to higher dimensional type
statements and to a study of special classes of vector bundles. As a quite surprising example,
one obtains statements of a very classical flavor about the ideals of special 
subvarieties in Jacobians. Some of these applications are explained in the present 
paper, while others, together with the necessary theoretical refinements, make the 
subject of the sequels \cite{pp1}, \cite{pp2}.

Let $X$ be an abelian variety of dimension $g$ over an algebraically closed field. 
Denote by $\hat{X}$ its dual, and let $\mathcal{P}$ be a Poincar\'e line 
bundle on $X\times \hat{X}$, properly normalized. The Fourier-Mukai functor is the 
derived functor associated to the functor which takes a coherent sheaf $\F$ on $X$ 
to $\hat{\mathcal{S}}(\F):=
{p_{\hat{X}}}_*(p_X^*\F\otimes \mathcal{P})$ on $\hat{X}$. We define Mukai 
regularity to be a condition on the cohomologies of the derived complex, weaker 
than the usual Index Theorem or Weak Index Theorem of \cite{mukai}. Concretely, a sheaf $\F$ on $X$ will 
be called Mukai regular, or simply \emph{$M$-regular}, if 
$${\rm codim}({\rm Supp}~R^i
\hat{\mathcal{S}}(\F))>i , ~\forall~ i>0.$$ 
Recall that a sheaf $\F$ is said to satisfy the Index Theorem (I.T.) with index $0$ if 
$$h^i(\F\otimes \alpha)=0, ~\forall~\alpha\in \Pic0, ~\forall~i>0,$$
thus sheaves satisfying I.T. with index $0$ are trivially $M$-regular. 
The main theme in what follows is that $M$-regularity has significant geometric 
consequences, and can be verified in 
practice in a variety of situations. As in many other settings, geometric information
is obtained via the global generation of suitably chosen sheaves, 
as we will see repeatedly below. 
This principle takes various forms, 
and the following is the main and most commonly used in this paper.

\begin{introtheorem}{\bf ($M$-regularity criterion)}
Let $\F$ be a coherent sheaf and $L$ an invertible sheaf supported on a 
subvariety $Y$ of the abelian variety $X$ (possibly $X$ itself). If both $\F$ and $L$ are 
$M$-regular as sheaves on $X$, then $\F\otimes L$ is globally generated.
\end{introtheorem}

By fixing a polarization $\Theta$ on $X$, one obtains a corollary which itself
provides a wide range of applications when combined with various vanishing theorems.
It is a generalization of a result on vector bundles proved in \cite{pareschi}
Theorem 2.1. Loosely speaking, it is precisely this extension to arbitrary coherent 
sheaves which allows one to attack a number of problems related to ideal sheaves 
or sheaves supported on subvarieties of $X$. 

\begin{introcorollary}{\bf (Theta-regularity criterion)}
Let $\F$ be a coherent sheaf on the polarized abelian variety $(X, \Theta)$. 
If $\F$ is $\Theta$-regular, i.e. $\F(-\Theta)$ is $M$-regular,
then $\F$ is globally generated. In particular, if $Y$ is a subvariety of $X$ such 
that $\I_Y(m\Theta)$ is $\Theta$-regular, then $\I_Y$ is cut out by divisors linearly 
equivalent to $m\Theta$.
\end{introcorollary}

\noindent
Based on a more refined notion of continuous global generation, we will see in fact 
that in the above the stronger fact that $\I_Y$ is cut out by divisors algebraically 
equivalent to $(m-1)\Theta$ also holds.
The Corollary above, and other Theta-regularity results,
are collected in \S6 in Theorem \ref{acm}, which can be seen as an abelian 
analogue of the well-known Castelnuovo-Mumford Lemma.

The $M$-regularity criterion is in turn a consequence of the following statement, which 
is the main technical result of the paper. The theorem is an extended generalization of 
a result on vector bundles appearing in various forms in work of Mumford, Kempf and Lazarsfeld
(cf. \cite{kempf1} or \cite{kempf2}), and explained in detail in \cite{pareschi}. 
The proof makes essential use of the full extent of Mukai's Fourier transform 
methods in a derived category setting. 

\begin{introtheorem}
Let $\F$ and $H$ be sheaves on $X$ such that $\F$ is $M$-regular and 
$H$ is locally free satisfying
I.T. with index $0$. Then, for any non-empty Zariski open set $U\subset \hat X$, the map
$${\mathcal M}_U:\bigoplus_{\xi\in U}H^0(X,\F\otimes P_\xi)\otimes H^0(X,H\otimes P_\xi^\vee)
\buildrel{\oplus m_\xi}\over
\longrightarrow H^0(X,\F\otimes H)$$
is surjective, where $m_{\xi}$ denote the multiplication maps on global sections.
\end{introtheorem}

Turning to applications, it is interesting to note in the first place that the  
apparently involved machinery of $M$-regularity has in fact consequences of a 
very elementary nature. The main application of this kind in the present 
paper is to the cohomological 
properties and the equations defining the images $W_d$ of symmetric products
of a curve $C$ in its polarized Jacobian $(J(C), \Theta)$, via Abel-Jacobi mappings.
The intuitive reason for the need of such methods is that naturally defined sheaves,
like the restriction $\OO_{W_d}(\Theta)$, are Mukai-regular, while not having 
better cohomological properties (cf. \S4).

\begin{introtheorem}
For any $1\leq d\leq g-1$, the twisted ideal sheaf is $3$-$\Theta$-regular, more 
precisely $\I_{W_d}(2\Theta)$ satisfies I.T. with index $0$.
\end{introtheorem}

\begin{introcorollary}
For any $1\leq d\leq g-1$, the ideal $\I_{W_d}$ is cut out by divisors
algebraically equivalent to $2\Theta$. It is also cut out by divisors 
linearly equivalent to $3\Theta$.
\end{introcorollary}

\noindent
As far as we are aware, the statement of the Corollary seems to be new even in the case $i=1$, 
i.e. for an Abel-Jacobi embedding of $C$ in its Jacobian. 
Another immediate consequence is an affirmative answer to Conjecture 10.3 in
\cite{oxbury}, on the sections of $2\Theta$ on symmetric products.

\begin{introcorollary}{\bf(Oxbury-Pauly conjecture)}
In the notation above, the following are true:
\newline
\noindent
(1) $h^0(C^d, u_d^*\OO(2\Theta))= \Sigma_{i=0}^d {g \choose i}$
\newline
\noindent
(2) $u_d^*:H^0(J(C),\OO(2\Theta))\rightarrow H^0(C^d, u_d^*\OO(2\Theta))$ is surjective.
\end{introcorollary}

In a rather different direction, in combination with Kawamata-Viehweg, Koll\'ar or 
Nadel vanishing type theorems (cf. \cite{rob}), our regularity criteria produce applications
with a higher dimensional flavor. In this paper we take only first steps in this 
direction, by looking at subvarieties of abelian varieties and their finite covers (but 
cf. also Remark \ref{maximal}),
while a more general study will be developed in \cite{pp2}. We collect 
in the following theorem some of the consequences obtained here (cf. \S5  
for these and further results).

\begin{introtheorem}
(i){\bf (Finite Albanese maps.)} Let $Y$ be a smooth irregular variety 
whose Albanese map is finite onto its 
image and let $L$ be a big and nef line bundle on $Y$. Then $\omega_Y^{\otimes 2}
\otimes L^{\otimes 2}$ is globally generated. In particular, if $Y$ is in addition a minimal 
variety of general type, then $\omega_Y^{\otimes 4}$ is globally generated. Moreover, 
if $\Theta$ is a polarization on ${\rm Alb}(Y)$, then $\omega_Y\otimes L\otimes
{\rm alb}_Y^*\OO(\Theta)$ is globally generated.
\newline
\noindent
(ii){\bf (Generalized Lefschetz Theorem.)} Under the same hypotheses, 
$\omega_Y^{\otimes 3}\otimes L^{\otimes 3}$ is very ample.
In particular, if $Y$ is also a minimal variety 
variety of general type, then $\omega_Y^{\otimes 6}$ is very ample.
\newline
\noindent
(iii){\bf (Relative base point freeness.)} Let $f:Y\rightarrow X$ be a proper, 
surjective morphism from a smooth 
variety $Y$ to an abelian variety $X$, and denote by $\omega_Y$ the canonical 
line bundle on $Y$. Then $f_* \omega_Y\otimes \OO_X(2\Theta)$ is globally generated,
for any ample divisor $\Theta$ on $X$.
\newline
\noindent
(iv){\bf (Mumford Lemma for multiplier ideals.)} Let $D$ be an effective $\QQ$-divisor
on a polarized abelian variety $(X, \Theta)$ 
and let $\J(D)$ be its associated multiplier ideal sheaf.
If $L$ is a divisor on $X$ such that $L-D$ is ample, then the sheaf
$\OO_X(\Theta + L)\otimes \J(D)$ is globally generated.
\end{introtheorem}
\noindent
The reader not familiar with the topic will find an explanation of the significance 
of (iv) above in \S6. Part (i) and (ii) are based on the fact that $\omega_Y\otimes L$
has an $M$-regular (in fact I.T.0) push-forward to the Albanese variety, which will also 
imply that $H^0(\omega_Y\otimes L)\neq 0$. This answers positively a nonvanishing conjecture 
of Kawamata \cite{kawamata2} in the case of varieties of this kind (cf. Corollary
\ref{nonvanishing}).

As one of the main applications of this theory, the concept of $M$-regularity 
allows for a unified point of view in the study of linear series on abelian 
varieties and their subvarieties. 
Here we discuss first cases of this assertion -- 
for example it is interesting to note that the classical Lefschetz and Ohbuchi 
very ampleness theorems and the fact that line 
bundles of degree $2g+1$ or more on curves are very ample are simply instances of the 
same phenomenon of producing $M$-regularity from numerical hypotheses (cf. \S4).
This can now be extended to a variety of numerical statements on 
arbitrary varieties which can be embedded in abelian varieties (for instance on $W_d$'s, 
to refer to the most direct generalization of the case of curves). 

Restricting the present discussion to linear series on abelian varieties, based on
the methods of this paper and those of \cite{pareschi}, in the 
sequel \cite{pp1} we will show among other things the following, improving on the 
well-known Lefschetz-Ohbuchi-Mumford-Kempf-Koizumi type results:

\begin{introtheorem}
(i) If $L$ is an ample line bundle on $X$ with $M$-regularity index $m(L)$, then
$L^{\otimes(k+2-m(L))}$ is $k$-jet ample.
\newline
\noindent
(ii) If $L$ is an ample line bundle with no base divisor on $X$, then the ideal of 
$X$ is cut out by quadratic equations in the embedding given by $L^{\otimes k}$, $k\geq 3$.
\newline
\noindent
(iii) Under the same hypothesis, the ideal of $X$ is cut out by quadratic and cubic 
equations in the embedding given by $L^{\otimes 2}$ (note that this is an embedding
by a theorem of Ohbuchi).
\end{introtheorem}
\noindent
The $M$-regularity index is a new invariant introduced in \cite{pp1}; we will not 
give the definition here.
Note only that if $L$ has no base 
divisor, then $m(L)\geq 1$. Statements (ii) and (iii) can be extended to arbitrary 
higher syzygies and also related to the invariant $m(L)$, improving on 
a conjecture of Lazarsfeld proved in \cite{pareschi}. 
In fact it will be seen in \cite{pp1} that essentially all the results on embeddings, 
higher order properties 
and defining equations for multiples of ample line bundles existing in the literature,
plus new results like the theorem above, can be subsumed in a more general 
framework depending on these invariants defined via $M$-regularity.  

Finally, it is important to emphasize that some of the underlying ideas in this work 
have separately existed in one form or another for quite some time. All our constructions 
rely heavily on Mukai's remarkable theory of the Fourier functor. On the other 
hand, it was Kempf who first realized that results in the style of Theorem 
\ref{amultiplication} provide a theta-group-free approach to statements on linear
series, based on vanishing theorems for vector bundles. The new ingredients 
brought in by the present approach are the inclusion of arbitrary coherent sheaves
in the general study, the systematic use of cohomological criteria for the global generation
of such sheaves, and the relaxation of strong vanishing conditions to the 
weaker $M$-regularity, all of which largely extend the range of applications.
It is here that this paper claims its main originality.

The paper is structured as follows. In Section 2 we introduce the basic 
definitions and prove the general 
$M$-regularity and multiplication criteria. Section 3 is devoted to a series 
of examples and first applications. Sections 4 and 
5 contain the main applications of the paper. In Section 4 we study the 
cohomological properties
and the equations defining the $W_d$'s in Jacobians, including the proof of the 
Oxbury-Pauly conjecture. In Section 5 we give various effective results for 
linear series on subvarieties of abelian varieties and more general irregular 
varieties.\footnote{We have recently found out about 
the paper \cite{chen} which has, among other things, an overlap with some of 
the results of this section.} 
Finally, in Section 6 we restrict to a discussion of Theta
regularity, drawing on a comparison with the usual Castelnuovo-Mumford regularity.

\noindent
\textbf{Acknowledgements.} Special thanks are due to R. Lazarsfeld, whose suggestions
at various stages have improved both authors' way of thinking about abelian varieties. 
We are also indebted to I. Dolgachev and B. van Geemen for some useful conversations, 
to C. Pauly for drawing our attention to the conjecture in \cite{oxbury}, and to 
J.A. Chen and C.D. Hacon for showing us their preprint \cite{chen} and for a remark 
improving the statement of Theorem \ref{finite_alb}(ii). 
The second author would like to thank B. van Geemen, K. O'Grady and University of 
Rome I for hospitality and financial support during a visit when this project was started.
Finally we thank the referee for careful suggestions which have allowed
us to fix some mistakes and substantially improve the exposition.

\noindent
\textbf{Background and notation.}
In what follows, unless otherwise specified, $X$ will be an abelian variety 
of dimension $g$ over an algebraically closed field.
We will always denote arbitrary coherent sheaves on $X$ with script letters (e.g. $\F$) and 
locally free sheaves with straight letters (e.g. $F$). Given a line bundle $L$, we 
write $B(L)$ for the base locus of the linear series $|L|$. 

We denote by $\hat{X}$ the dual abelian variety, 
which will often be identified with $\Pic0$.
Given a point $\xi\in \hat{X}$, by $P_\xi\in \Pic0$ we understand the corresponding 
line bundle via this identification. By
$\mathcal{P}$ we denote a Poincar\'e
line bundle on $X\times \hat{X}$, normalized such that 
$\mathcal{P}_{|X\times\{0\}}$
and $\mathcal{P}_{| \{0\}\times \widehat{X}}$ are trivial.

We briefly recall the Fourier-Mukai setting, referring to \cite{mukai} for details.
To any coherent sheaf $\F$ 
on $X$ we can associate the sheaf ${p_{2}}_{*}({p_{1}}^{*}\F\otimes\mathcal{P})$ on
$\widehat{X}$, where $p_1$ and $p_2$ are the natural projections on $X$ and $\hat{X}$. 
This correspondence gives a functor
$\hat{\mathcal{S}}: {\rm Coh(X)}\rightarrow {\rm Coh}(\widehat{X})$. 
If we denote by ${\bf D}(X)$ and 
${\bf D}(\hat{X})$ the derived categories of Coh$(X)$ and Coh$(\hat{X})$, then 
the derived functor $\mathbf{R}\hat{\mathcal{S}}:{\bf D}(X)\rightarrow {\bf D}(\hat{X})$ 
is defined and called the \emph{Fourier functor},
and one can consider $\mathbf{R}\mathcal{S}:{\bf D}(\hat{X})\rightarrow {\bf D}(X)$ 
in a similar way. Mukai's main result \cite{mukai} Theorem 2.2 
is that $\mathbf{R}\hat{\mathcal{S}}$ establishes 
an equivalence of categories between ${\bf D}(X)$ and ${\bf D}(\hat{X})$.

Let $R^j\hat{\mathcal{S}}(\F)$ be the cohomologies of the derived complex 
$\mathbf{R}\hat{\mathcal{S}}(\F)$. 
A coherent sheaf $\F$ on $X$ satisfies W.I.T. (the \emph{weak index theorem}) with index $i$ if 
$R^j\hat{\mathcal{S}}(\F)=0$ for all $j\neq i$. It satisfies the stronger I.T. (the
\emph{index theorem}) with index $i$ if 
$H^{i}(\F\otimes\alpha)=0$ for all $\alpha\in {\rm Pic}^{0}(X)$ and all $i\neq j$.
By the base change theorem, in this situation $R^{j}\hat{\mathcal{S}}(\F)$ is locally free. 
If $\F$ satisfies W.I.T. with index $i$, $R^{i}\hat{\mathcal{S}}(\F)$ is denoted by 
$\hat{\F}$ and called the $\it{Fourier~
transform}$ of $\F$. Note that then $\mathbf{R}\hat{\mathcal{S}}(\F)\cong \hat{\F}[-i]$.

\section{\bf Mukai regularity, global generation and continuous global generation}

\noindent
Let $X$ be an abelian variety of dimension $g$.
Given a coherent sheaf $\F$ on $X$, we denote 
$$S^i(\F):={\rm Supp}(R^i\hat{\mathcal S}(\F)).$$

\begin{definition}
A coherent sheaf $\F$ on $X$ is called \emph{Mukai-regular} (or simply 
\emph{$M$-regular}) if
${\rm codim}(S^i(\F))>i$ for any $i=1,\ldots ,g$ (where, for $i=g$, 
this means that $S^g(\F)$ is empty).
\end{definition}

\begin{example}
Sheaves satisfying I.T. with index $0$ are the simplest 
examples of $M$-regular sheaves. A first important class is that of ample line bundles 
on $X$.
\end{example}

\begin{remark}
From the definition of the Fourier transform we see that there is always an inclusion
$S^i(\F)\subset V^i(\F)$, where $V^i(\F)$ is the \emph{cohomological 
support locus} (cf. 
\cite{green}):
$$V^i(\F):= \{~\xi ~|~h^i(\F\otimes P_{\xi})\neq 0\}\subset {\rm Pic}^0(X).$$
Consequently, $M$-regularity is achieved if in particular
$${\rm codim}(V^i(\F))>i {\rm ~for~ any~} i=1,\ldots ,g.$$
It will be this property that we will usually check in applications.
\end{remark}

The key point in what follows is that $M$-regularity is a cohomological condition
which has significant geometric consequences via the global generation of suitable sheaves.
The main result is the following:

\begin{theorem}\label{areg}
Let $\F$ be a coherent sheaf and $L$ an invertible sheaf supported on a 
subvariety $Y$ of the abelian variety $X$ (possibly $X$ itself). If both $\F$ and $L$ are 
$M$-regular as sheaves on $X$, then $\F\otimes L$ is globally generated.
\end{theorem}

The proof of Theorem \ref{areg} essentially occupies the rest of this section; more precisely, 
it is a formal consequence of Proposition \ref{product}
and Proposition \ref{cont_glob_gen} below. 
The key technical theorem, which will also find independent use in later sections, 
generalizes in several ways a result of 
Mumford-Kempf-Lazarsfeld explained in \cite{pareschi} \S2.

\begin{theorem}\label{amultiplication}
Let $\F$ and $H$ be sheaves on $X$ such that $\F$ is $M$-regular and 
$H$ is locally free satisfying
IT with index $0$. Then, for any non-empty Zariski open set $U\subset \hat X$, the map
$${\mathcal M}_U:\bigoplus_{\xi\in U}H^0(X,\F\otimes P_\xi)\otimes H^0(X,H\otimes P_\xi^\vee)
\buildrel{\oplus m_\xi}\over
\longrightarrow H^0(X,\F\otimes H)$$
is surjective, where $m_{\xi}$ denote the multiplication maps on global sections.
\end{theorem}

\begin{proof}
Note that when there is no danger of confusion, we will avoid writing explicitly the 
names of the maps involved.
The conclusion of the theorem is equivalent to the injectivity of the dual map:
$$H^0(\F\otimes H)^\vee\longrightarrow \prod_{\xi\in U}H^0(\F\otimes P_\xi)^\vee\otimes 
H^0(H\otimes P_\xi^\vee)^\vee,$$
or, by Serre duality (and the fact that $H$ is locally free), with the injectivity 
of the co-multiplication map:
\begin{equation}\label{first_form}
{\rm Ext}^g(\F,H^\vee)\longrightarrow\prod_{\xi\in U}{\rm Hom}
(H^0(\F\otimes P_\xi),H^g(H^\vee\otimes P_\xi)).
\end{equation}

We will concentrate on showing (\ref{first_form}). 
The proof requires the language of derived categories.
Let us first reindex, for technical reasons, $R^j\hat{\mathcal{S}}(\F)$
by $R^{-j}\hat{\mathcal{S}}(\F)$, $-g\leq j\leq 0$. Note also that since $H$ 
satisfies I.T. with index $0$, by Serre duality $H^\vee$ 
satisfies I.T with index $g$, and so $\widehat{H^\vee}=R^g\hat{\mathcal{S}}(H^{\vee})$
is locally free.
We will freely identify $\widehat{H^\vee}$ with the one-term complex 
${\bf R}\hat{\mathcal{S}}(H^{\vee})[g]$.

\noindent
\emph{Claim 1.} {\it There exists a bounded fourth quadrant cohomological spectral sequence 
$$E^{ij}_2 ={\rm Ext}^i(R^{-j}\hat{\mathcal{S}}(F),\widehat{H^\vee})\Rightarrow
H^{i+j}={\rm Ext}^{i+j}_{{\bf D}(\hat X)}({\bf R}\hat{\mathcal{S}}(F),\widehat{H^\vee}).$$}

\noindent
Consider the functors $F={\rm Hom}(\circ , \widehat{H^\vee})$ and $G= \hat{\mathcal{S}}$.
We first apply the derived functor ${\bf R}G$ to $\F$ to obtain an object in 
${\bf D}(\hat X)$. We can then consider a standard  hypercohomology 
spectral sequence (cf. \cite{weibel} 5.7.9 and 10.8.3): 
$$E^{ij}_2 = (R^i F)(H^j {\bf R}G (\F))\Rightarrow {\mathbb R}^{i+j}F({\bf R}G (\F)),$$
where ${\mathbb R}^i F$ are the right hyper-derived functors of $F$ (in this case 
the hyperexts). 
This is precisely the sequence in Claim 1. Note that we have reindexed 
the $R^j\hat{\mathcal{S}}(\F)$
in order to agree with the usual sign convention used for the Hom cochain complex 
\cite{weibel} 2.7.4.. The rest of the assertions are clear.

\noindent
\emph{Claim 2.} {\it The spectral sequence in Claim 1 induces a natural inclusion
$$H^{0+0}= {\rm Hom}_{{\bf D}(\hat X)}({\bf R}\hat{\mathcal{S}}(\F), \widehat{H^\vee}) 
\hookrightarrow E_{2}^{00}= {\rm Hom}(R^0\hat{\mathcal{S}}(\F),\widehat{H^\vee}).$$}

\noindent
Note that the $M$-regularity assumption on $\F$ guarantees that 
${\rm Ext}^i(R^{-j}\hat{\mathcal{S}}(\F),\widehat{H^\vee})=0$ for $i \leq -j$ and $i>0$, 
so that the only $E^{i,j}_\infty$ term such that $i+j=0$ is $E_\infty^{00}$.
This implies that
$$H^{0+0}={\rm Hom}_{{\bf D}(\hat X)}({\bf R}\hat{\mathcal{S}}(\F), 
\widehat{H^\vee})\cong E_\infty^{00}.$$
But since we have a fourth quadrant spectral sequence, the differentials 
coming into $E^{00}_p$ are always zero, so we easily get a chain of inclusions
$$E_{\infty}^{00}\subset \ldots\subset E_{3}^{00}\subset E_{2}^{00}.$$

\noindent
\emph{Claim 3.} {\it There is a natural map of $\OO_{\hat X}$-modules
$$\phi:{\rm Ext}^g(\F,H^\vee)\otimes\OO_{\hat X}\longrightarrow 
\mathcal{H}om(R^0\hat{\mathcal{S}}(\F), \widehat{H^\vee})$$
such that, for $\xi$ in a suitable open set $V\subset \hat X$, the induced maps on the fibers
$$\phi(\xi):{\rm Ext}^g(\F,H^\vee)\longrightarrow \mathcal{H}om(R^0\hat{\mathcal{S}}(\F),
\widehat{H^\vee})(\xi)$$
are the co-multiplication maps in (\ref{first_form}).}

\noindent
By Mukai's duality result \cite{mukai} Theorem 2.2, we have
$${\rm Ext}^g (\F, H^\vee) \cong {\rm Hom}_{{\bf D}(X)}(\F,H^\vee [g])
\cong {\rm Hom}_{{\bf D}(\hat X)}({\bf R}\hat{\mathcal{S}}(\F), 
{\bf R}\hat{\mathcal{S}}(H^\vee)[g])$$
$$\cong {\rm Hom}_{{\bf D}(\hat X)}({\bf R}\hat{\mathcal{S}}(\F), \widehat{H^\vee}[g-i(H^\vee)]) 
\cong {\rm Hom}_{{\bf D}(\hat X)}({\bf R}\hat{\mathcal{S}}(\F), \widehat{H^\vee}).$$
We thus get a map at the level of global sections
\begin{equation}\label{global_scns}
\Phi=H^0(\phi): {\rm Ext}^g(\F,H^\vee)\longrightarrow {\rm Hom}(R^0\hat 
{\mathcal{S}}(\F),\widehat{H^\vee}),
\end{equation}
induced by the spectral sequence map in Claim 2. This can be extended in turn to a sheaf map as in 
the statement of the claim (since the first sheaf is the trivial bundle). 
It is natural, and not hard to check, that this 
map coincides with the above co-multiplication maps at the general point.

Finally, $\mathcal{H}om(R^0\hat{\mathcal{S}}(\F),R^g\hat{\mathcal{S}}(H^\vee))$
is a torsion-free sheaf (note that 
$R^g\hat{\mathcal{S}}(H^\vee)$ is locally free), so it is a standard fact that 
the injectivity of the map in (\ref{first_form}) is equivalent
to the injectivity of the map $\phi$ of Claim 3 at the $H^0$ level, in other 
words that of the map $\Phi$ in (\ref{global_scns}).
This is precisely the statement of Claim 2, and the proof is completed.
\end{proof}

\begin{remark}
If in the theorem above we impose the stronger condition that $\F$ satisfy I.T. with index $0$, 
then the argument above shows the degeneration of the spectral sequence at the $E_2$ 
level, which implies that the map $\Phi$ is in fact an isomorphism.  
\end{remark}

\begin{remark}
For the purposes of the present paper we only need Theorem \ref{amultiplication}
in the case when $H$ is an ample line bundle on $X$. There are however interesting 
applications of the full statement, for example to the study of semihomogeneous vector 
bundles, which we will describe in \cite{pp2}.
\end{remark}

\begin{corollary}\label{finite}
Let $\F$ and $H$ as in the previous theorem. Then there is a positive integer $N$ 
such that for any general
$\xi_1,\ldots,\xi_N\in \hat X$, the map
$$\bigoplus_{i=1}^NH^0(X,\F\otimes P_\xi)\otimes H^0(X,H\otimes P_\xi^\vee)
\buildrel{\oplus m_{\xi_i}}\over\longrightarrow H^0(X,\F\otimes H)$$
is surjective.
\end{corollary}
\begin{proof}
Let $U$ be an open set
of $\hat X$ on which the rank of the map $m_\xi$ is constant. The surjectivity of the map
${\mathcal M}_U$ implies that the family of linear subspaces $\{{\rm Im}(m_\xi)\}_{\xi\in U}$ 
linearly spans
the vector space $H^0(\F\otimes H)$. But then there is an $N$ such that $N$ general such 
subspaces span $H^0(\F\otimes H)$.
\end{proof}

In the same spirit we have a ``preservation of vanishing'' statement which will be used 
in the subsequent sections to deduce vanishing results from $M$-regularity results.

\begin{proposition}\label{freg_vanishing}
Let $\F$ be an $M$-regular coherent sheaf on $X$ and $H$ a locally free sheaf satisfying I.T. with 
index $0$. Then $\F\otimes H$ satisfies I.T. with index $0$.
\end{proposition} 
\begin{proof}
For an arbitrary $\alpha\in \Pic0$, since $H$ is locally free we have 
$$H^i(\F\otimes H\otimes \alpha) \cong {\rm Ext}^i((H\otimes \alpha)^{\vee}, \F).$$
Applying again Mukai's duality theorem \cite{mukai} Theorem 2.2  we get
$${\rm Ext}^i((H\otimes \alpha)^{\vee}, \F) = {\rm Hom}_{{\bf D}(X)}((H\otimes \alpha)^{\vee}, \F[i])\cong
{\rm Hom}_{{\bf D}(\hat{X})}({\bf R}\hat{\mathcal{S}}((H\otimes \alpha)^{\vee}), {\bf R}\hat{\mathcal{S}}(\F)[i])$$
$$\cong {\rm Hom}_{{\bf D}(\hat{X})}({((H\otimes \alpha)^{\vee})}^{\widehat{}}, {\bf R}\hat{\mathcal{S}}(\F)[i+g])
\cong {\rm Ext}^{i+g}_{{\bf D}(\hat{X})}(((H\otimes \alpha)^{\vee})^{\widehat{}}, {\bf R}\hat{\mathcal{S}}(\F)).$$
As in the proof of Theorem \ref{amultiplication}, there is a bounded cohomological spectral sequence
$$E^{ij}_2 ={\rm Ext}^i(((H\otimes \alpha)^{\vee})^{\widehat{}}, R^{j}\hat{\mathcal{S}}(F))\Rightarrow
H^{i+j}={\rm Ext}^{i+j}_{{\bf D}(\hat X)}(((H\otimes \alpha)^{\vee})^{\widehat{}}, {\bf R}\hat{\mathcal{S}}(F)).$$
Note now that the $M$-regularity assumption on $X$ implies that 
$${\rm Ext}^i(((H\otimes \alpha)^{\vee})^{\widehat{}}, R^{j}\hat{\mathcal{S}}(F))=0 {\rm ~for~} i+j>g,$$ 
so in particular $E^{i,j}_\infty = 0$ in the same range. This immediately implies that 
$${\rm Ext}^{i+g}_{{\bf D}(\hat{X})}(((H\otimes \alpha)^{\vee})^{\widehat{}}, {\bf R}\hat{\mathcal{S}}(\F))
= 0 {\rm ~for~} i>0$$
as claimed.
\end{proof}

\medskip
The fact that Theorem \ref{amultiplication} produces geometric
statements is best explained 
by introducing the intermediate notion of continuous global generation, which can 
be in fact defined on an arbitrary irregular 
variety. As it will be clear from the discussion below, 
a good understanding of this concept provides in turn 
global generation statements via Proposition \ref{product}.

\begin{definition}
Let $Y$ be an irregular variety. We define a sheaf $\F$ on $Y$ to 
be \emph{continuously globally generated} if for any non-empty open subset 
$U\subset {\rm Pic}^{0}(Y)$ the sum of evaluation maps 
$$\bigoplus_{\alpha\in U} H^0(\F\otimes \alpha)\otimes \alpha^\vee \longrightarrow \F$$ 
is surjective. 
\end{definition}

\begin{remark}
As in Corollary \ref{finite}, this is equivalent to the existence of a number $N$ 
such that the finite sum of multiplication maps taken over $N$ general line bundles 
in ${\rm Pic}^{0}(Y)$ is surjective.
The intuitive significance of this notion comes from the fact that it models 
precisely the behavior of line bundles on abelian varieties: the "continous
system" of all divisors algebraically equivalent to a given one on an abelian variety
is "base point free", i.e. the intersection of all divisors in the class is empty. 
Thus continuous global generation is a generalization to arbitrary coherent 
sheaves on irregular varieties. Moreover, for line bundles it means precisely the 
same as in the case of abelian varieties: the intersection of the divisors in 
$|L\otimes \alpha|$ is empty as $\alpha$ varies over an open subset of ${\rm Pic}^0(Y)$.
\end{remark}

\medskip

\begin{proposition}\label{product}
Let $Y$ be a subvariety of an irregular variety $X$, $\F$ a coherent sheaf and $L$ 
a line bundle supported on  
$Y$, both continuously globally generated as sheaves on $X$. 
Then $\F\otimes L$ is globally generated.
\end{proposition}
\begin{proof}
We denote the inclusion of $Y$ in $X$ by $i$, but for simplicity we use
the same notation $\F$ or $L$ even when we consider them as torsion 
sheaves on $X$. We then have $H^0(X, \F\otimes P_{\xi})\cong 
H^0(Y, \F\otimes i^*P_{\xi})$ and the analogous statement for $L$. 
Consider the following commutative diagram, obtained by 
alternating the order of the obvious evaluation and multiplication maps:
$$\xymatrix{
\bigoplus_{i=1}^NH^0(\F\otimes P_{\xi_i})\otimes H^0(L\otimes 
P_{\xi_i}^\vee)\otimes \OO_X \ar[r] \ar[d] &
H^0(\F\otimes L)\otimes\OO_X \ar[d]^{ev_{\F\otimes L}} \\
\bigoplus_{i=1}^NH^0(\F\otimes P_{\xi_i})\otimes L\otimes 
P_{\xi_i}^\vee \ar[r] & \F\otimes L } $$
where $N$ is chosen so that for general $\xi_1,\ldots, \xi_N\in {\rm Pic}^0(X)$ 
the bottom horizontal
map is onto. It follows that the support of Coker($ev_{\F\otimes L}$) is contained 
in the union of the base loci $B(L\otimes i^*P_{\xi_i}^\vee)$ on $Y$, for $\xi_1,\ldots,
\xi_N\in {\rm Pic}^0(X)$ general, or equivalently the loci where $L\otimes P_{\xi}^\vee$ 
are not generated by global sections on $X$ individually, 
since $L$ has rank at most $1$ everywhere. 
But the continuous global generation condition implies precisely 
that the intersection of these loci is empty (cf. also the Remark above), which easily 
gives that the intersection of such unions is empty.
\end{proof}

Finally the key point, using the new language, 
is that $M$-regular sheaves are always continuously globally generated.
More precisely, we have the following:

\begin{proposition}\label{cont_glob_gen}{\bf ($M$-regularity implies 
continuous global generation.)}
If $\F$ is $M$-regular, then 
there is a positive integer $N$ such that for general $\xi_1,\ldots ,\xi_N\in
\hat X$, the sum of twisted evaluation maps 
$$\bigoplus_{i=1}^N H^0(\F\otimes P_{\xi_i})\otimes P_{\xi_i}^\vee \longrightarrow \F$$ 
is surjective.
\end{proposition}
\begin{proof}
We apply Theorem \ref{amultiplication} with $H$ a line bundle,  
sufficiently ample so that $\F\otimes H$ is globally generated. 
As in the proof of Proposition \ref{product}, consider the commutative diagram:
$$\xymatrix{
\bigoplus_{i=1}^NH^0(\F\otimes P_{\xi_i})\otimes H^0(H\otimes 
P_{\xi_i}^\vee)\otimes \OO_X \ar[r] \ar[d] &
H^0(\F\otimes H)\otimes\OO_X \ar[d] \\
\bigoplus_{i=1}^NH^0(\F\otimes P_{\xi_i})\otimes H \otimes 
P_{\xi_i}^\vee \ar[r] & \F\otimes H } $$
By Corollary \ref{finite} and the choice of $H$, the top-right composition 
must be surjective, so the same must hold for the bottom map, which is what we 
wanted to show. 
\end{proof}

\section{\bf First examples and applications of $M$-regularity}

\begin{example}\label{freg_curves}
A line bundle $L$ on a smooth curve $C$ of genus $g\geq 1$ is $M$-regular (via an 
Abel-Jacobi embedding $C\subset J(C)$) if and only if $d={\rm deg}L \geq g$.
Indeed, the $M$-regularity of $L$ is equivalent by Riemann-Roch to the fact that
$W_d^{d-g+1}$ has codimension at least $2$ in ${\rm Pic}^d(C)$. But this is easily 
seen to be equivalent to $d\geq g$.
\end{example}

\begin{example}
More generally, let $W_d \subset J(C)$ be an Abel-Jacobi embedding given by 
the choice of a line bundle of degree $d$ on $C$, and consider the natural map 
$\pi_d: C^d\rightarrow W_d$. Let also $N_1,\ldots N_d$ 
be line bundles of degree $g$ on $C$, and on $W_d$ consider the sheaf 
$$\F:={\pi_d}_*(N_1\boxtimes\ldots \boxtimes N_d).$$
Then one can see similarly that $\F$ is $M$-regular as a sheaf on $J(C)$.
From Theorem \ref{areg} and Example \ref{higher_d} below, one gets then 
immediately the following 

\begin{corollary}
The sheaf $\F\otimes \OO_{W_d}(\Theta)$ is globally generated on $W_d$.
\end{corollary}

\noindent
This example is of interest in the sense that it is a priori not clear (at least 
to us) how to prove this global generation by standard methods. Note that 
$\F\otimes \OO_{W_d}(\Theta)$ is the push-forward of the line bundle 
$(N_1\boxtimes\ldots \boxtimes N_d)\otimes \pi_d^*\OO(\Theta)$ on $C_d$, which is 
of the form 
$$[A\otimes N_1)\boxtimes \ldots \boxtimes (A
\otimes N_d)]\otimes \OO(-\Delta),$$
where $A$ is a line bundle of degree $g+d-1$ on $C$ and $\Delta$ is the union of 
the diagonal divisors (cf. \S5).
\end{example}

\begin{example}\label{higher_d}
If $\Theta$ is a theta divisor on the Jacobian of $C$, then $\OO_{W_d}(\Theta)$
is $M$-regular. This fact will be important in \S5. The proof is not immediate,
but see Proposition \ref{freg_wd}.
\end{example}

\noindent
\textbf{Linear series -- an introduction to \cite{pp1}.}
In a different direction, our $M$-regularity and global generation results
provide a unified approach to very ampleness (in a first stage) statements on 
abelian varieties and their subvarieties. We will need the following standard lemma.

\begin{lemma}\label{very_ampleness}
Let $L$ be a line bundle on a variety $Y$. Then $L$ is very ample if and 
only if $L\otimes \I_y$ is globally generated for all $y\in Y$.
\end{lemma}

\noindent
As an immediate consequence of Theorem \ref{areg} and Lemma \ref{very_ampleness}, we get:

\begin{corollary}\label{va_subvar}
Let $Y$ be a subvariety of an abelian variety $X$. Let $L$ and $M$ be line bundles 
on $Y$ such that $L$ is $M$-regular and $M\otimes \I_y$ is $M$-regular for all $y\in Y$, 
both as sheaves on $X$. Then $L\otimes M$ is very ample.
\end{corollary}

\begin{lemma}\label{base_divisor}
If $L$ is an ample line bundle with no base divisor on an abelian variety $X$, 
then $L\otimes \I_x$ is $M$-regular for any $x\in X$.
\end{lemma}
\begin{proof}
The cohomological support locus $V^1(L\otimes \I_x)$ is 
the locus of $\xi\in \hat{X}$
corresponding to $P_\xi \in \Pic0$ such that $h^1(L\otimes P_\xi
\otimes \I_x)\neq 0$. Since $L\otimes P_\xi \otimes \I_x\cong t^*_{\nu}L\otimes 
\I_x\cong t^*_{\nu}(L\otimes \I_{x-\nu})$ for 
some $\nu\in X$, and since $h^1(L)=0$, this is precisely isomorphic 
to the locus where $L$ fails to be globally generated.
But now $L$ has no base divisor, so this implies that 
$V^1(L\otimes \I_x)$ has codimension at least $2$. Also $V^i(L\otimes \I_x)$ 
are obviously empty for $i\geq 2$, so $L\otimes \I_x$ is $M$-regular.
\end{proof}

\begin{example}
To give a flavor of the range of applications that can be derived as  
consequences of Corollary \ref{va_subvar}, we explain how some of the most basic 
facts on curves and abelian varieties are incorporated in this theory. 
For the abelian varieties statements below cf. \cite{lange} Ch.4 \S5.

\begin{corollary}
Let $X$ be an abelian variety  and $C$ a smooth projective curve. Then:
\newline
\noindent
(i)(\emph{Lefschetz Theorem}) If $L$ is an ample line bundle
on $X$, then $L^{\otimes 3}$ is very ample.
\newline
\noindent
(ii)(\emph{Ohbuchi's Theorem}) If $L$ is an ample line bundle on $X$ 
with no base divisor, then $L^{\otimes 2}$ is very ample.
\newline
\noindent
(iii) If $L$ is a line bundle of degree at least $2g+1$ on $C$, then 
$L$ is very ample.
\end{corollary}
\begin{proof}
By Theorem \ref{va_subvar}, (i) and (ii) follow if we prove that
$L^{\otimes 2}\otimes \I_x$, respectively $L\otimes \I_x$, are $M$-regular. 
But this is a consequence of Lemma \ref{base_divisor}, since in the first case 
$L^{\otimes 2}$ is globally generated and in the second case $L$ has no base divisor 
by assumption. For part (iii) note that a line bundle of degree $2g+1$ on $C$ can 
be written as $A\otimes B$, where $A$ has degree $g$ and $B$ has degree $\geq g+1$.
By Example \ref{freg_curves}, $A$ and $B(-x)$ are $M$-regular for all $x\in X$, and 
so again by Theorem \ref{areg} this implies that $A\otimes B(-x)$ is globally generated.
\end{proof}
\end{example}

\begin{remark}
The above corollary is just intended to be an amusing illustration of the fact that 
well-known and seemingly unrelated results are in fact realizations of the same principle.
On the other hand, and more importantly, the proof is a toy version of a general approach, 
based on the notion of $M$-regularity, which provides
a unified treatment of essentially all the known geometric and syzygetic statements on 
multiples of ample line bundles on abelian varieties, and also produces new 
basic statements on higher order properties of such embeddings, projective normality 
and defining equations, as explained in the introduction.
This is the main topic of the sequel \cite{pp1} to this paper. 
\end{remark}

\section{\bf The equations defining the $W_d$'s in Jacobians and a conjecture 
of Oxbury and Pauly}

Let $C$ be a smooth curve of genus 
$g\geq 3$, and denote by $J(C)$ the Jacobian of $C$. Let $\Theta$ be a theta
divisor on $J(C)$, and $C_d$ the $d$-th
symmetric product of $C$, for $1\leq d\leq g-1$. Consider 
$$u_d :C_d\longrightarrow J(C)$$ 
to be an Abel-Jacobi mapping of the symmetric product (depending on the choice 
of a line bundle of degree $d$ on $C$), and denote by $W_d$ 
the image of $u_d$ in $J(C)$.
In this section we show how the methods developed in this paper allow one to
understand the "regularity" of the ideal sheaf $\I_{W_d}$ (cf. \S6 for the precise 
meaning of this). This gives a 
positive answer to a conjecture of Oxbury and Pauly \cite{oxbury}
(cf. Corollary \ref{op_conj} below) and bounds in turn the degrees of the 
Theta-equations defining $W_d$ in $J(C)$.

\begin{theorem}\label{w_d}
For any $1\leq d\leq g-1$, the twisted ideal sheaf 
$\I_{W_d}(2\Theta)$ satisfies I.T. with index $0$.
\end{theorem}

A first application of this result is to defining equations.
One would like to know what multiples of $\Theta$ cut out the ideal of 
$W_d$ in $J(C)$, in other words for what $m$ the sheaf $\I_{W_d}\otimes \OO_J(m\Theta)$
is globally generated (or better, continuously globally generated). 
The answer to this question turns out to be independent of $d$, namely:

\begin{corollary}\label{cubic_theta}
For any $1\leq d\leq g-1$, the ideal $\I_{W_d}$ is cut out by divisors
algebraically equivalent to $2\Theta$. 
Moreover, $\I_{W_d}$ is also cut out by divisors linearly equivalent to $3\Theta$.
\end{corollary}
\begin{proof}
The first statement means precisely that $\I_{W_d}(2\Theta)$ 
is continuously globally generated. But this follows immediately from Proposition
\ref{cont_glob_gen}, since by Theorem \ref{w_d} $\I_{W_d}(2\Theta)$ satisfies
I.T. with index $0$, hence it is $M$-regular.
On the other hand, by Theorem \ref{areg}
this implies that $\I_{W_d}(3\Theta)$ is globally generated, which gives
the second statement.
\end{proof}

As mentioned above, Theorem \ref{w_d} also gives an affirmative answer to Conjecture 10.3 in
\cite{oxbury}, stated there in the slightly more restrictive case when $\Theta=\Theta_{\kappa}$,
where $\kappa$ is a theta-characteristic on $C$.

\begin{corollary}\label{op_conj}{\bf(Oxbury-Pauly conjecture)}
In the notation above the following are true:
\newline
\noindent
(1) $h^0(C^d, u_d^*\OO(2\Theta))= \Sigma_{i=0}^d {g \choose i}$
\newline
\noindent
(2) $u_d^*:H^0(J(C),\OO(2\Theta))\rightarrow H^0(C^d, u_d^*\OO(2\Theta))$ is surjective.
\end{corollary}

\begin{proof}
Both statements follow immediately from the vanishing of $h^i(\I_{W_d}\otimes \OO(2\Theta))$ 
for all $i>0$, combined with the computation of $\chi(u_d^*\OO(2\Theta))$ given in 
\cite{oxbury} Proposition 10.1(3). But this vanishing (together with that for any other 
translate) means precisely that $\I_{W_d}(2\Theta)$  
satisfies I.T. with index $0$, as stated in Theorem \ref{w_d}.
\end{proof}

\bigskip
\noindent
A key step in the proof of Theorem \ref{w_d} is the following 
(cf. also Example \ref{higher_d}):

\begin{proposition}\label{freg_wd}
(i) The sheaf $\OO_{W_d}(\Theta)$ is $M$-regular on $J(C)$.
\newline 
(ii) We have $h^0(W_d, \OO_{W_d}(\Theta)\otimes P_{\xi})=1$ for 
$\xi\in {\rm Pic}^0(J(C))$ generic.
\end{proposition}
\begin{proof}
We will prove the two statements at the same time, by computing explicitly 
the corresponding cohomology groups. 
The Abel-Jacobi map has projective spaces as fibers,
so the vanishings we are interested in follow as soon as 
we prove them after pulling back via $u_d$.
We are then interested in computing the cohomology groups 
$$H^i(C_d, u_d^*\OO(\Theta)\otimes u_d^*P_{\xi}), ~\forall~\xi\in {\rm Pic}^0(J(C)).$$
For this we appeal to the technical results proved in \cite{izadi} Appendix 3.1, 
which we use freely below. If we denote by $\pi_d: C^d \rightarrow J(C)$ the corresponding 
desymmetrized Abel-Jacobi map, we have:
$$\pi_d^*\OO(\Theta)\cong (A\boxtimes\ldots\boxtimes 
A)\otimes \OO(-\Delta),$$
where $\Delta$ is the union of all the diagonal divisors in $C^d$ and $A$ is a line 
bundle of degree $g+d-1$. We then have the identification:
$$H^i(u_d^*\OO(\Theta)\otimes u_d^*P_{\xi})
\cong S^i H^1(C, A\otimes \xi)  \otimes \wedge^{d-i}
H^0(C, A\otimes \xi),$$
obtained as the skew-symmetric part of the cohomology group 
$H^i((A\otimes \xi)^{\boxtimes d})$ under the action of $S_d$.
Now we only have to compute the loci on which these cohomology groups do not
vanish, for each $0 \leq i\leq d$. Since $\xi$ is a general line bundle on $C$ 
of degree $0$, $A\otimes \xi$ is a general line bundle
of degree $g+d-1$. Thus statement (ii) follows immediately, since 
$H^1(C, A\otimes \xi)=0$ and so 
$$H^0(C, A\otimes \xi)=d.$$
On the other hand, for $i>0$, if $H^i(u_d^*\OO(\Theta)\otimes u_d^*P_{\xi})$ is nonzero 
then $h^1(A\otimes\xi)=h^0(\omega_C\otimes A^{-1}\otimes 
\xi^{-1})\neq 0$.
This locus is isomorphic to $W_{g-d-1}$, and so it is of codimension $d+1$,
which gives (i).
\end{proof}

\begin{proof}(\emph{of Theorem \ref{w_d}})
From the standard cohomology sequence associated to the exact sequence
$$0\longrightarrow \I_{W_d}(2\Theta)\longrightarrow \OO_{J(C)}(2\Theta)
\longrightarrow \OO_{W_d}(2\Theta)\longrightarrow 0,$$
we see that the theorem follows as soon as we prove the following two statements, 
where this time for simplicity we denote by $\Theta$ any of its translates:
\begin{itemize}
\item $H^i(\OO_{W_d}(2\Theta))=0, ~\forall~i>0$.
\item The restriction map $H^0(\OO_{J(C)}(2\Theta))
\rightarrow H^0(\OO_{W_d}(2\Theta))$ is surjective.
\end{itemize}
The first statement means precisely that $\OO_{W_d}(2\Theta)$ satisfies I.T. with 
index $0$, which follows immediately from Proposition \ref{freg_wd}(i) and 
Proposition \ref{freg_vanishing}.

For the second statement we will appeal directly to Theorem \ref{amultiplication}. 
Note that for any open subset $U\in \widehat{J(C)}$ we have a commutative 
diagram as follows, where the vertical maps are the natural restrictions. 
$$\xymatrix{
\bigoplus_{\xi \in U} H^0(\OO_{J(C)}(\Theta_\xi))\otimes 
H^0(\OO_{J(C)}(\Theta_{-\xi})) \ar[r] \ar[d] &
H^0(\OO_{J(C)}(2\Theta)) \ar[d] \\
\bigoplus_{\xi \in U} H^0(\OO_{W_d}(\Theta_\xi))\otimes 
H^0(\OO_{J(C)}(\Theta_{-\xi})) \ar[r] & H^0(\OO_{W_d}(2\Theta)) } $$
Now Theorem \ref{amultiplication} says that the bottom horizontal map 
is surjective. On the other hand, by Proposition \ref{freg_wd}(ii) we can choose 
the open set $U$ such that the left vertical map is an isomorphism. (Note that we cannot 
include all the translates $\Theta_\xi$ such that 
$h^0(\OO_{W_d}(\Theta_\xi))=1$, but certainly those that satisfy this property plus 
the fact that $W_d$ is not contained in $\Theta_\xi$.) This in turn implies that the 
right vertical map is surjective, which is precisely our statement.
\end{proof}

\section{\bf {Applications via vanishing theorems}}

In combination with various vanishing theorems, like those of Kawamata-Viehweg, 
Koll\'ar or Nadel, the $M$-regularity criterion 
produces effective geometric results on linear series  and special coherent 
sheaves on subvarieties of
abelian varieties and more general irregular varieties, in the spirit of 
\cite{kollar2}.
We assume here that we are working over a field of characteristic zero.
In all the applications below, $M$-regularity will be satisfied in the strong 
form of I.T. with index $0$. 

\noindent
\textbf{Effective global generation and very ampleness.} 
The following result treats the case of (finite covers of) 
subvarieties of abelian varieties (but see also Remark \ref{maximal} below).

\begin{theorem}\label{finite_alb}
(i) Let $Y$ be a smooth irregular variety whose Albanese map is finite onto its 
image and let $L$ be a big and nef line bundle on $Y$. Then $\omega_Y^{\otimes 2}
\otimes L^{\otimes 2}$ is globally generated. In particular, if $Y$ is in addition a minimal 
variety of general type, then $\omega_Y^{\otimes 4}$ is globally generated.
\newline
\noindent
(ii) Under the same hypotheses,   
$\omega_Y^{\otimes 3}\otimes L^{\otimes 3}$ is very ample.
In particular, if $Y$ is also a minimal variety 
variety of general type, then $\omega_Y^{\otimes 6}$ is very ample.
\end{theorem}
\begin{proof}
(i) By Proposition \ref{product}, the result follows if we prove that $\omega_Y\otimes L$ is 
continuously globally generated. As with global generation, it is easy to see
that the continuous global generation of $\omega_Y\otimes L$ is implied by that
of ${{\rm alb}_Y}_*(\omega_Y\otimes L)$, since ${\rm alb}_Y$ is finite. The Kawamata-Viehweg
vanishing theorem on $Y$, and again the fact that ${\rm alb}_Y$ is finite, imply that
$$h^i({{\rm alb}_Y}_*(\omega_Y\otimes L)\otimes \alpha)=
h^i(\omega_Y\otimes L\otimes \alpha)= 0, ~\forall~ i>0, ~\forall~\alpha\in {\rm Pic}^0(Y)
={\rm Pic}^0({\rm Alb}(Y)).$$
This means that ${{\rm alb}_Y}_*(\omega_Y\otimes L)$ satisfies I.T. with index $0$, hence 
it is continuously globally generated by Proposition \ref{cont_glob_gen}.
\newline
\noindent 
(ii) By Lemma \ref{very_ampleness}, Proposition \ref{product} and the fact that 
$\omega_Y\otimes L$ is continuously globally generated, 
it would  be enough to show that 
$\omega_Y^{\otimes 2}\otimes L^{\otimes 2}\otimes \I_y$ is continuously globally 
generated for all $y\in Y$. As above, this is implied by the same statement for the 
push-forward to ${\rm Alb}(Y)$, which in turn follows as long as
${{\rm alb}_Y}_*(\omega_Y^{\otimes 2}\otimes L^{\otimes 2}\otimes \I_y)$ satisfies 
I.T. with index $0$. But this is immediate, since as $\omega_Y\otimes L$ is 
continuously globally generated, it is nef and so $\omega_Y\otimes L^{\otimes 2}$ 
is big and nef. Consequently
$$H^i(\omega_Y^{\otimes 2}\otimes L^{\otimes 2}\otimes \I_y \otimes \alpha)=0$$
(for all $\alpha\in {\rm Pic}^{0}(Y)$), trivially for $i\geq 2$, and 
because of the global generation of (deformations of) 
$\omega_Y^{\otimes 2}\otimes L^{\otimes 2}$ for $i=1$.
\end{proof}

\begin{remark}
The result above can be stated more generally for Gorenstein varieties with
rational singularities, on which the vanishing theorem still applies.
\end{remark}

\begin{remark}\label{maximal}
Although we have preferred to state the result above in the most compact 
form, it is important to mention that the same statments hold (and with essentially 
the same argument) in the more naturally occuring situation of 
varieties of maximal Albanese dimension, at least outside the 
exceptional loci, i.e. outside the higher dimensional fibers of the generically 
finite Albanese map. 
\end{remark}

The point in the above is that given a big and nef line bundle $L$ on an
irregular variety $Y$ with finite Albanese map, $\omega_Y\otimes L$ is continuously 
globally generated. This implies in particular both that it is nef and that 
it has global sections. Consequently this gives a positive answer to Kawamata's 
nonvanishing conjecture \cite{kawamata2} on varieties of this kind
(note that the nefness hypothesis in the conjecture is superfluous in this case).

\begin{corollary}\label{nonvanishing}
Let $Y$ be a smooth irregular variety with finite Albanese map and let $L$ be a 
big and nef line bundle on $Y$. Then $H^0(\omega_Y\otimes L)\neq 0$.
\end{corollary}

The proof of Theorem \ref{finite_alb} produces also a similar result
depending this time on a fixed arbitrary polarization on the target 
abelian variety.

\begin{proposition}\label{finite_theta_alb}
Let $f:Y\rightarrow X$ be a finite morphism from a projective variety to a polarized
abelian variety $(X, \Theta)$, and $\F$ a coherent sheaf on $Y$. If $f_{*}\F$ is 
$M$-regular, then $\F\otimes f^*\OO_X(\Theta)$ is globally generated. In particular, 
if $Y$ is smooth and 
$L$ is big and nef on $Y$, then 
$\omega_Y\otimes L\otimes f^*\OO_X(\Theta)$
is globally generated and
$\omega_Y^{\otimes 2}\otimes L^{\otimes 2}\otimes f^*\OO_X(\Theta)$ is very ample, 
where $\Theta$ is an arbitrary polarization on $X$.
\end{proposition}

On a slightly different note,
Fujita's problem on the global generation of adjoint bundles has a
relative version, which we learned of from \cite{kawamata}. Since on abelian varieties
the square of an ample line bundle is always globally generated, it is 
natural to ask whether there is an analogous uniform statement in the 
relative setting. This again follows quickly from Theorem \ref{areg}.

\begin{proposition}\label{kaw}
Let $f:Y\rightarrow X$ be a proper, surjective morphism from a smooth 
variety $Y$ to an abelian variety $X$, and denote by $\omega_Y$ the canonical 
line bundle on $Y$. Then $f_* \omega_Y\otimes \OO_X(2\Theta)$ is globally generated,
for any ample divisor $\Theta$ on $X$.
\end{proposition}

\begin{proof}
This follows from Theorem \ref{areg} via a special case of Koll\'ar's
theorem on direct images of dualizing sheaves \cite{kollar}, which states that
$$h^i(f_*\omega_Y\otimes L) = 0,~\forall i>0,$$
for any ample line bundle $L$ on $X$.
\end{proof}

\noindent
\textbf{An optimal Mumford type Lemma for multiplier ideals.}
Vanishing theorems for coherent 
sheaves often lead to global generation statements via the Castelnuovo-Mumford
Lemma (cf. \S6). In particular, a general such statement for 
multiplier ideal sheaves is getting to be known as Mumford's Lemma. This has found 
some very important applications, for example in Siu's proof of the deformation invariance
of plurigenera, or to Fujita type statements on the volume of big divisors
(cf. \cite{rob} \S4.3 for a general discussion). We show here that on abelian varieties one can 
do substantially better: the analogy to keep in mind is the difference 
between the general Fujita type statements and the best possible global generation 
bounds on abelian varieties (like the global generation of 
$\OO_X(2\Theta)$). For the general theory of multiplier ideals we refer the 
reader to \cite{rob}.

Let $D$ be an effective $\QQ$-divisor on $X$ and denote by 
$\J(D)$ the associated multiplier ideal sheaf.

\begin{proposition}\label{mult_ideals}
If $L$ is a divisor on $X$ such that $L-D$ is ample, then the sheaf
$\OO_X(\Theta + L)\otimes \J(D)$ is globally generated.
\end{proposition}
\begin{proof}
By Theorem \ref{areg}, it is sufficient to have the 
vanishing:
$$h^i(\OO_X(L)\otimes\J(D)\otimes \alpha)=0, ~\forall~ i>0, ~\forall~ \alpha\in {\rm Pic}^0(X).$$
Now by assumption $\OO_X(L-D)\otimes \alpha$ is ample (and of course $\omega_X = \OO_X$),
so this follows from the Nadel vanishing theorem (cf. \cite{rob} Theorem 4.6).
\end{proof}

Although not in direct relationship with our criterion, in order to complete the 
picture we note also that the abelian varieties variant of the nonvanishing theorem
for multiplier ideals \cite{rob} \S4.3 holds in a very strong form.

\begin{proposition}
Under the assumptions of Proposition \ref{mult_ideals}
$$H^0(\OO_X(L)\otimes \J(D))\neq 0.$$
\end{proposition}
\begin{proof}
Since $\OO_X(L-D)\otimes \alpha$ is ample for all $\alpha\in {\rm Pic}^0(X)$, the Nadel 
vanishing theorem implies that $\OO_X(L)\otimes \J(D)$ satisfies I.T. with index $0$.
By Corollary \ref{cont_glob_gen} it is then continuously globally generated, so 
it has sections.
\end{proof}

One can formulate analogous statements for multiplier ideals of $\QQ$-divisors
on smooth subvarieties of abelian varieties. One such version, whose proof is 
essentially that of Proposition \ref{mult_ideals} together with Proposition
\ref{product}, can be stated as follows:

\begin{proposition}
Let $Y$ be a smooth subvariety of $X$ and $D$ an effective $\QQ$-divisor on $Y$.
Assume that $L$ is a divisor such that $L-D$ is big and nef, and $M$ is 
a continuously globally generated divisor on $Y$. Then the sheaf
$$\OO_Y(K_Y + L + M)\otimes \J(D)$$ 
is globally generated.
\end{proposition}

\section{\bf Theta regularity and an ``abelian'' Castelnuovo-Mumford Lemma}

The notion of $M$-regularity might not be an a priori obvious analogue of the usual
notion of regularity on projective spaces. 
There is however a particular instance of the general theory, depending on 
a fixed polarization $\Theta$, 
which is more intuitive and is already responsible for some of the 
applications above.
We will call this \emph{Theta regularity} in what follows.
The point of this section is to see that the theory of Theta regularity 
is strikingly similar to the theory of 
Castelnuovo-Mumford regularity for coherent sheaves on (subvarieties of) projective 
spaces. 

\begin{definition}
A coherent sheaf $\F$ on a polarized abelian variety $(X, \Theta)$ is 
called $m$-$\Theta$-\emph{regular} 
if $\F((m-1)\Theta)$ is $M$-regular.
If $\F$ is $0$-$\Theta$-regular, we will simply call it $\Theta$-\emph{regular}.
\end{definition}

\begin{example}
Many of the results of the previous two sections can be phrased in the language 
of Theta regularity. For example, $\OO_{W_d}$ is $2$-$\Theta$-regular, 
while Theorem \ref{w_d} implies that $\I_{W_d}$ is $3$-$\Theta$-regular.
\end{example}

We collect in the following theorem
a number of results which show the similarity mentioned above. This can be 
considered a Theta-regularity version of the Castelnuovo-Mumford Lemma. For the 
usual Castelnuovo-Mumford Lemma cf. e.g. \cite{mumford3}, Lecture 14.

\begin{theorem}\label{acm}
Let $\F$ be a $\Theta$-regular coherent sheaf on $X$. Then:
\newline
\noindent
(1) $\F$ is globally generated.
\newline
\noindent
(2) $\F$ is $m$-$\Theta$-regular for any $m\geq 1$.
\newline
\noindent
(3) The multiplication map
$$H^0(\F(\Theta))\otimes H^0(\OO(k\Theta))\longrightarrow H^0(\F((k+1)\Theta))$$
is surjective for any $k\geq 2$.
\end{theorem}

\begin{proof}
Part (1) is a particular case of Theorem \ref{areg}, while 
(2) follows immediately from Proposition \ref{freg_vanishing}, as $\OO(\Theta)$
satisfies I.T. with index $0$.
For (3) note that, for any open subset $U\subset \hat{X}$, we have a commutative diagram
$$\xymatrix{ 
\bigoplus_{\xi\in U}H^0(\F(-\Theta)\otimes P_\xi)\otimes H^0({\OO}_X(2\Theta)
\otimes P_\xi^\vee)\otimes H^0({\OO}_X(k\Theta)) \ar[r] \ar[d] &
H^0(\F(\Theta))\otimes H^0({\OO}_X(k\Theta)) \ar[d] \\
\bigoplus_{\xi\in U}H^0(\F(-\Theta)\otimes
P_\xi)\otimes H^0(\OO_X((k+2)\Theta)\otimes P_\xi^\vee) \ar[r] & 
H^0(\F((k+1)\Theta)) }$$
Now the left vertical map is surjective if $U$ is chosen such that 
Ohbuchi's projective normality theorem (cf. \cite{ohbuchi}) applies 
for any $\xi\in U$. In addition, since $\F(-\Theta)$ 
is assumed to be $M$-regular, the bottom horizontal map is also surjective by 
Theorem \ref{amultiplication}. This implies that the right vertical map is surjective.
\end{proof}

\begin{remark}
The numerical analogy with Castelnuovo-Mumford regularity
is not perfect in part (3) of Theorem \ref{acm}. 
It is easy to see though that the statement in (3) is optimal, as it follows 
for example by considering $\F$ equal to $\OO(2\Theta)$ (when we cannot make 
$k=1$). This particular chioce 
also shows that this "abelian" Castelnuovo-Mumford statement contains the 
well-known Mumford type 
projective normality results for multiples of ample line bundles (cf. \cite{mumford2}).
Note also that in this form part (3) generalizes to Theta-regular coherent sheaves
a weaker statement about locally free sheaves satisfying I.T. with index $0$, 
which follows as a consequence of \cite{pareschi} Theorem 3.8.  
\end{remark}

\noindent
\textbf{Bounds on $\Theta$-regularity in terms of defining equations.}
The Castelnuovo-Mumford regularity of a subvariety in projective space is a measure 
of the complexity of the computations involving its ideal (cf. e.g. \cite{bayer}).
It is natural to ask for bounds on this
invariant in terms of the degrees of defining equations for the subvariety. This has 
been optimally achieved for smooth subvarieties of projective spaces in \cite{bertram}.
Completely similar arguments (which we will not repeat here) can be used to give a bound 
on the $\Theta$-regularity of a smooth subvariety of an abelian variety in  
terms of the degrees of its pluritheta generators. This can be quite naturally seen
as a converse of the cohomological global generation criterion Theorem \ref{acm} (1) in the 
case of ideal sheaves of smooth subvarieties.

\begin{theorem}\label{converse}
Let $Y$ be a codimension $e$ smooth subvariety of the polarized abelian variety 
$(X,\Theta)$, and 
assume that the ideal of $Y$ in $X$ is cut out by $d\Theta$-equations (i.e. the 
sheaf $\I_Y\otimes \OO_X(d\Theta)$ is globally generated). Then $\I_Y$ is 
$ed$-$\Theta$-regular, more precisely $\I_Y\otimes \OO_X((ed+1)\Theta)$ satisfies
I.T. with index $0$. Moreover, the $\Theta$-regularity of $\I_Y$ is precisely 
$ed$ if and only if $Y$ is a complete intersection of $d\Theta$-hypersurfaces in $X$.
\end{theorem}

\noindent 
Theorem \ref{converse} also has a more general version, where the degrees of the defining 
$\Theta$-equations are allowed to vary, along the lines of \cite{bertram} Corollary 4.

\providecommand{\bysame}{\leavevmode\hbox to3em{\hrulefill}\thinspace}

\end{document}